\documentclass[11pt]{amsart}
\usepackage{amssymb,amsmath,amsfonts,amscd,euscript}

\newcommand{\nc}{\newcommand}

\nc{\slt}{\mathfrak{sl}_2}
\nc{\s}{\mathfrak{sl}}
\nc{\suth}{\widehat{\mathfrak{su}}(2)}
\nc{\gl}{\mathfrak{gl}}
\nc{\GL}{\mathfrak{GL}}
\nc{\g}{\mathfrak{g}}
\nc{\gh}{\widehat{\mathfrak{g}}}
\nc{\h}{\mathfrak{h}}
\nc{\hh}{\widehat{\mathfrak{h}}}
\nc{\la}{\lambda}
\nc{\slth}{\widehat{\slt}}
\nc{\C}{\mathbb C }
\nc{\R}{\mathbb R }
\nc{\Z}{\mathbb Z }
\nc{\N}{\mathbb N }
\nc{\al}{\alpha}
\nc{\be}{\beta}
\nc{\ve}{\varepsilon}
\nc{\ch}{{\mathop {\rm ch}}}
\nc{\Id}{{\mathop {\rm Id}}}
\nc{\Tr}{{\mathop {\rm Tr}\,}}
\nc{\U}{{\mathop {\rm U}}}
\nc{\bra}{\langle}
\nc{\ket}{\rangle}
\nc{\ld}{\ldots}
\nc{\cd}{\cdots}
\nc{\hk}{\hookrightarrow}
\nc{\n}{\mathfrak{n}}
\nc{\un}{\mathfrak{u}}
\nc{\bo}{\mathfrak{b}}
\nc{\T}{\otimes}
\nc{\wt}{\widetilde}
\nc{\bn}{{\bf n}}
\nc{\veps}{\varepsilon}

\nc{\qb}[2]{\genfrac{(}{)}{0pt}{}{#1}{#2}_q}
\nc{\at}[2]{\genfrac{}{}{0pt}{}{#1}{#2}}

\nc{\beq}{\begin{equation}}
\nc{\eeq}{\end{equation}}

\newtheorem{theo}{Theorem}[section]
\newtheorem*{theo*}{Theorem}
\newtheorem{lem}{Lemma}[section]
\newtheorem{prop}{Proposition}[section]
\newtheorem{cor}{Corollary}[section]
\newtheorem{rem}{Remark}[section]
\newtheorem{conj}{Conjecture}[section]
\newtheorem{defin}{Definition}[section]

\begin{document}
\author{E.Feigin}
\title
[The PBW filtration]
{The PBW filtration}

\address{Evgeny Feigin:\newline
{\it Tamm Theory Division, Lebedev Physics Institute,
Russian Academy of Sciences,\newline 
Russia, 119991, Moscow, Leninski prospect, 53}\newline 
and \newline
{\it Independent University of Moscow,\newline 
Russia, Moscow, 119002, Bol'shoi Vlas'evskii, 11}}
\email{evgfeig@gmail.com}

\begin{abstract}
In this paper we study the PBW filtration on irreducible integrable highest 
weight representations of affine Kac-Moody algebra $\gh$. 
The $n$-th space of this
filtration is spanned with the vectors $x_1\dots x_s v$, where  $x_i\in\gh$, $s\le n$  and $v$ is 
a highest weight vector. For the vacuum module we give a 
conjectural description of the corresponding adjoint graded space  in
terms of generators and relations. For $\g$ of the type $A_1$ we prove our conjecture
and derive the fermionic formula for the graded character.
\end{abstract}

\maketitle

\section*{Introduction}
Let $\g$ be a Kac-Moody Lie algebra of finite or affine type and $\U(\g)$ be the 
universal enveloping algebra of $\g$. For the dominant integral
weight $\la$ let $L_\la$ be the corresponding irreducible highest weight representation
with highest weight $\la$ and highest weight vector $v_\la$. We consider a 
filtration $\U(\g)_s$ on $\U(\g)$ defined by 
$$
\U(\g)_0=\C 1,\ \U(\g)_{s+1}=\U(\g)_s + \mathrm{span}\{gu:\ g\in \g, u\in \U(\g)_s\}.
$$  
This filtration induces a filtration $F_s=\U(\g)_s\cdot v_\la$ on $L_\la$. Define an
associated graded space 
$$
L_\la^{gr}=F_0\oplus F_1/F_0\oplus F_2/F_1\oplus\dots
$$
We define the graded  $u$-character by
$$
\ch_u L_\la^{gr}=\sum_{s\ge 0} u^s \ch (F_s/F_{s-1}),
$$
where we set $F_{-1}=0$ and $\ch$ denotes the usual character with respect to the Cartan 
subalgebra of $\g$. For example for $\g$ of the type $A_1$ and its arbitrary finite-dimensional
representation all the spaces $F_s/F_{s-1}$ are one-dimensional. The PBW-filtration 
for $\g$ of the type $A_1^{(1)}$ and its level $1$ representations was used in
\cite{FFJMT} for the study of $\phi_{1,3}$-field in Virasoro minimal theories.  

Let us briefly describe our approach to the study of the spaces $L_\la^{gr}$. Let
$\n_-\hk\g$ be the nilpotent subalgebra of creating operators, i.e. 
\beq
L_\la=\U(\n_-)\cdot v_\la.
\eeq
The action of $\n_-$ on $L_\la$ induces the action of the abelian algebra 
$\n_-^{ab}$ on $L_\la^{gr}$ ($\n_-^{ab}$ is isomorphic to $\n_-$ as a vector space).
Therefore the space $L_\la^{gr}$ is generated from the vector $v_\la$ with 
the action
of $\U(\n_-^{ab})$, which is isomorphic to the polynomial algebra. This allows
to describe $L_\la^{gr}$ as a quotient of the polynomial algebra on $\n_-$ by some
ideal.

In this paper we study the case of an affine Kac-Moody algebra $\gh$ and its vacuum 
representations. 
Let $\g$ be a finite-dimensional simple Lie algebra and $\gh$ be the corresponding
affine algebra, i.e. the central extension of $\g\T\C[t,t^{-1}]$.
We fix  its vacuum level $k$ 
representation $L_k$ with a highest weight vector $v_k$. In this case one has
$
(\g\T \C[t])\cdot v_k=0
$
and therefore the space $L_k^{gr}$ is generated from the highest weight vector 
with an action of the universal enveloping algebra of the abelian algebra 
$\g^{ab}\T t^{-1}\C[t^{-1}]$. This gives
\beq
L_k^{gr}\simeq  \U(\g^{ab}\T t^{-1}\C[t^{-1}])/I_k,
\eeq 
where $I_k$ is some ideal. We describe the ideal $I_k$ in the following way.
We note that $\g=\g\T 1$ annihilates the highest weight vector $v_k$. Therefore
the action of $\g$ on $L_k$ induces a structure of $\g$-module on $L_k^{gr}$ and $I_k$.
(This structure plays an important role in our approach, see the conjecture below.
For non vacuum modules $L_\la^{gr}$ is not $\g$-module anymore. This means that the 
ideal of relations
for general highest weight is to be described in other terms then in the vacuum case.)   
We give a conjectural description of $I_k$
below. The proof is given in the paper for the case $\g=\slt$.

The adjoint action of $\g$ on itself endows the space
$\U(\g^{ab}\T t^{-1}\C[t^{-1}])$ with a structure of $\g$-module.
Let $\theta$ be the longest root of $\g$, $f_\theta\in\g$ be an element of the weight
$\theta$ and
$$
f_\theta(z)=\sum_{n\ge 0} z^n (f_\theta\T t^{-n-1})
$$ 
be the corresponding current. We recall (see for example \cite{BF}) that the series
$f_\theta(z)^{k+1}$ acts by zero on $L_k$
(the coefficients of the $z$-expansion of $f_\theta(z)^{k+1}$ are acting by zero).
\begin{conj}
We have an equality
\beq
I_k=(\U(\g)\oplus \U(\g^{ab}\T t^{-1}\C[t^{-1}]))\cdot \mathrm{span}\{
\text{ coefficients of }  f_\theta(z)^{k+1} \},
\eeq
i.e. $I_k$ is the minimal $\U(\g)$-stable ideal which contains the coefficients of 
$f_\theta(z)^{k+1}$.
\end{conj}
We verify this conjecture for $\g=\slt$ (the level $1$ case was studied in \cite{FFJMT} by
different method). Combining the description
of $I_k$ with the vertex operator realization technique we obtain the fermionic formula 
for the $u$-character of $L_k^{gr}$ for $\g=\slt$.

The paper is organized in the following way:

In Section $1$ we fix our notations and recall some constructions from the representation
theory of affine Kac-Moody algebras and bosonic vertex operator algebras.

In Section $2$ we study the adjoint graded space $L_1^{gr}$.

In Section $3$ we describe the ideal of relations in $L_k^{gr}$ for the 
general level and obtain the fermionic formula for the graded character.

\noindent {\it Acknowledgments.} 
Research of EF was partially supported by the 
RFBR Grants 06-01-00037, 07-02-00799 and LSS 4401.2006.2. 
This paper was written during the authors stay at the Max Planck Institute
for Mathematics. EF thanks the Institute for hospitality.

\section{Preliminaries and definitions}
\subsection{Kac-Moody Lie algebras, representations and filtrations.}
Let $\g$ be a finite-dimensional simple Lie algebra with a Cartan decomposition
$\g=\n\oplus\h\oplus \n_-$.We consider the corresponding affine Kac-Moody algebra\
\beq
\gh=\g\T \C[t,t^{-1}]\oplus \C K\oplus\C d,
\eeq 
where $K$ is a central element and $d$ is a degree operator, $[d,x\T t^i]=-ix\T t^i$
for all $x\in\g$.  The affine algebra $\gh$ admits the Cartan decomposition
$\gh=\widehat{\n}\oplus\widehat{\h}\oplus\widehat{\n_-}$, where
\begin{gather*}
\widehat{\n}=\n\T 1\oplus \g\T t\C[t],\\ 
\widehat{\h}=\h\T 1\oplus \C K\oplus \C d,\\
\widehat\n_-=\n_-\T 1\oplus \g\T t^{-1}\C[t^{-1}].
\end{gather*}
Let $\la\in\widehat \h^*$ be an integral dominant weight and
$L_\la$ be the corresponding irreducible highest weight representation of $\gh$ with
highest weight $\la$ and highest weight vector $v_\la$. Then
$$
\widehat{\n} \cdot v_\la=0,\ h v_\la=\la(h) v_\la\ \forall h\in\widehat{\h},\ 
L_\la=\U(\widehat\n_-)\cdot v_\la,
$$
where $\U(\widehat\n_-)$ is a universal enveloping algebra. 
The number $\la(K)$ is called the level of $L_\la$. We let $L_k$ to
denote the level $k$ vacuum representation (the restriction of the highest 
weight of $L_k$ to $\h\T 1$ vanishes).

We define an increasing filtration $\U(\widehat\n_-)_s$ on the universal enveloping
algebra by the following rule:
$$
\U(\widehat\n_-)_0=\C \cdot 1,\ \U(\widehat\n_-)_{s+1}=\U(\widehat\n_-)_s +
\widehat\n_- \U(\widehat\n_-)_s.
$$
This filtration induces an increasing  filtration  $F_s$ on $L_\la$:
$$F_s=\U(\widehat\n_-)_s\cdot v_\la.$$ 
The filtration $F_\bullet$ is called the Poincare-Birkhoff-Witt filtration 
(the PBW filtration for short).

\begin{defin}
Define $L_\la^{gr}$ as an adjoint graded space with respect to the filtration
$F_s$, i.e. 
\beq
L_\la^{gr}=F_0\oplus\bigoplus_{s>0} F_s/F_{s-1}.
\eeq 
Define the $u$-character of $L_\la^{gr}$ as
\beq
\ch_u L_\la^{gr}=\sum_{s\ge 0} u^s \ch (F_s/F_{s-1}),
\eeq 
where we set $F_{-1}=0$ and $\ch$ denotes the usual character with respect to 
$\widehat \h^*$.
\end{defin}

\begin{rem}
\label{nab}
Let $\widehat \n_-^{ab}$ be an abelian Lie algebra which coincides with $\widehat \n_-$ as
a vector space. Then $L_\la^{gr}$ carries a natural structure of representation of 
$\widehat \n_-^{ab}$.
\end{rem}

\begin{rem}
The definition above is valid in the case of finite-dimensional algebras as well,
providing a filtration on finite-dimensional modules.
\end{rem} 

\begin{rem}
For the general weight $\la$ the spaces $F_s$ are not stable under the action of an algebra
$\g=\g\T 1\hk \gh$. This motivates the following modification of the filtration above.
Let 
$$
\widetilde F_0=\U(\n_-\T 1)\cdot v_\la,\ 
\widetilde F_{s+1}=\widetilde F_s + \widehat\n_- \widetilde F_s.
$$
Here $F_0=\C v_\la$ is replaced by $\widetilde F_0$, which is a finite-dimensional 
irreducible representation of $\g$.
We denote the corresponding
adjoint graded object by $\widetilde L_\la^{gr}$. It is obvious that for any $s$ the space
$\widetilde F_{s}$ is $\g$-invariant. We hope to return
to the study of this filtration elsewhere.   
\end{rem}

\subsection{The $\slt$ case}
Let $\g=\slt$ and let $e,h,f$ be its standard basis. 
Then $f$ spans $\n_-$, $h$ spans
$\h$ and $e$ spans $\n$. For $x\in\slt$ we set $x_i=x\T t^i$. 
Let $v_k$ be a highest weight vector of the vacuum level $k$ module $L_k$. 

We now recall a construction of the $ehf$-basis of $L_k$ from 
\cite{FKLMM}. 

\begin{defin}
\label{ehfdef}
A monomial of the form 
\beq
\label{form}
\dots f_{-n}^{a_n} h_{-n}^{b_n} e_{-n}^{c_n}\dots  
f_{-1}^{a_1} h_{-1}^{b_1} e_{-1}^{c_1}
\eeq 
is called an ordered monomial. In addition it is called a $ehf$-monomial 
if it satisfies the following conditions:
\begin{enumerate}
\item[(a)] $a_i+a_{i+1}+b_{i+1}\le k$ for $i> 0$, 
\label{a}
\item[(b)] $a_i+b_{i+1}+c_{i+1}\le k$ for $i> 0$, 
\label{b}
\item[(c)] $a_i+b_i+c_{i+1}\le k$ for $i> 0$,
\label{c}
\item[(d)] $b_i+c_i+c_{i+1}\le k$ for $i> 0$. 
\label{d}
\end{enumerate}
\end{defin}

\begin{theo}
\label{basis}
The set $\{m\cdot v_k\}$, where $m$  runs over the set of 
$ehf$-monomials provides a
basis of $L_k$.
\end{theo}
 
\begin{rem}
\label{pict}
The following picture from \cite{FKLMM} illustrates the set of $ehf$-monomials:

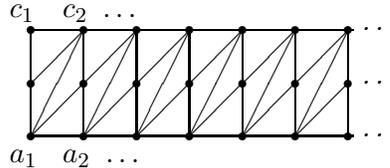
\begin{center}
\begin{picture}(180,70)

\multiput(40,20)(0,40){2}{\line(1,0){122}}
\multiput(40,20)(20,0){7}{\line(0,1){40}}
\multiput(40,20)(20,0){6}{\line(1,1){20}}
\multiput(40,20)(20,0){6}{\line(1,2){20}}
\multiput(40,40)(20,0){6}{\line(1,1){20}}
\multiput(165,20)(0,20){3}{\dots}
\multiput(40,20)(20,0){7}{\circle*{3}}
\multiput(40,40)(20,0){7}{\circle*{3}}
\multiput(40,60)(20,0){7}{\circle*{3}}
\put(32,10){$a_1$} 
\put(52,10){$a_2\ \dots$}
\put(32,65){$c_1$}
\put(52,65){$c_2 \ \dots$}

\end{picture}
\end{center}              

Namely one considers a set of monomials $(\ref{form})$ 
such that the sum of exponents over any
triangle (corresponding to the conditions (a)--(d))
is less than or equal to $k$.  

\end{rem}

\subsection{Lattice vertex operator algebras and affine Kac-Moody algebras}
In this section we recall main properties of lattice vertex
operator algebras (VOA for short) and their
principal subspaces.
The main references are \cite{K2}, \cite{BF},
\cite{dong}, \cite{FK}.

Let $Q$ be a lattice of finite rank equipped with a symmetric bilinear form
$(\cdot,\cdot): Q\times Q\to \Z$ such that $(\al,\al)>0$ for  all $\al\in 
Q\setminus \{0\}$. Let $\h=Q\T_{\Z} \C$. The form  $(\cdot,\cdot)$ induces
a bilinear form on $\h$, for which we use the same notation. Let 
$\h\T\C[t,t^{-1}]\oplus \C K$
be the corresponding multi-dimensional Heisenberg algebra with the bracket
$$[\al\T t^i, \be\T t^j]=i\delta_{i,-j}(\al,\be) K,\  
[K,\al\T t^i]=0,\ \al,\be\in\h.$$
For $\al\in\h$ define the 
Fock representation $\pi_\al$ of the Heisenberg algebra generated by a vector $|\al\ket$
 such that
$$(\be\T t^n) |\al\ket=0,\ n>0; \qquad (\be\T 1) |\al\ket=(\be,\al) |\al\ket;
\qquad
K|\al\ket=|\al\ket.$$

We now define a VOA $V_Q$ associated with $Q$. We deal only with an even
case, i.e.  $(\al,\al)\in 2\Z$ for all $\al\in Q$ 
(in the general case the construction leads to the 
so called super VOA). As a
vector space 
$$V_Q\simeq \bigoplus_{\al\in Q} \pi_\al.$$
The $q$-degree on $V_Q$ is defined by 
\begin{equation}
\label{defqdeg}
\deg_q |\al\ket= \frac{(\al,\al)}{2},\quad \deg_q (\be\T t^n)=-n.
\end{equation}
The main ingredients of the VOA structure 
on $V_Q$ are bosonic vertex 
operators
$\Gamma_\al(z)$ which correspond to highest weight vectors $|\al\ket$.
One sets
\begin{equation}
\label{defvo}
\Gamma_\al(z)=S_\al z^{\al\T 1} \exp(-\sum_{n<0} \frac{\al\T t^n}{n} z^{-n})
\exp(-\sum_{n>0} \frac{\al\T t^n}{n} z^{-n}),
\end{equation}
where $z^{\al\T 1}$ acts on $\pi_\be$ by $z^{(\al,\be)}$ and 
the operator $S_\al$ is defined by
$$S_\al |\be\ket=c_{\al,\be} |\al+\be\ket;\quad
[S_\al,\be\T t^n]=0, \ \al,\be\in\h,$$
where $c_{\al,\be}$ are some non vanishing constants.
The Fourier decomposition is given by
$$\Gamma_\al(z)=\sum_{n\in\Z} \Gamma_\al(n) z^{-n-(\al,\al)/2}.$$
In particular,
\begin{equation}
\label{htoh}
\Gamma_\al(-(\al,\al)/2-(\al,\be))|\be\ket=c_{\al,\be}|\al+\be\ket.
\end{equation}
One of the main properties of vertex operators is the following
commutation relation:
\begin{equation}
\label{vertop}
[\al\T t^n,\Gamma_\be(z)]=(\al,\be)z^n \Gamma_\be(z).
\end{equation}
Another important formula describes the product of two vertex operators
\begin{multline}
\Gamma_\al(z)\Gamma_\be(w)=(z-w)^{(\al,\be)} S_\al S_\be 
z^{(\al+\be)\T 1}\times\\ 
\exp(-(\sum_{n<0} \frac{\al\T t^n}{n} z^{-n} + \frac{\be\T t^n}{n} w^{-n}))
\exp(-(\sum_{n>0} \frac{\al\T t^n}{n} z^{-n} +  \frac{\be\T t^n}{n} w^{-n})).
\end{multline}
This leads to the proposition:

\begin{prop}
\label{vorel}
\begin{equation}
(\Gamma_\al(z))^{(k)}(\Gamma_\be(z))^{(l)}=0 \text{ if } k+l<(\al,\be),
\end{equation}
where the superscript $(k)$ denotes the $k$-th derivative of the corresponding
series. In addition if $(\al,\be)=0$ then 
$$\Gamma_\al(z)\Gamma_\be(z) \text{ is proportional to } \Gamma_{\al+\be}(z).$$
\end{prop}

We now recall the Frenkel-Kac construction which provides a vertex 
operator realization of the basic representation of the affine algebra
$\gh$ for $\g$ of the types $A$, $D$ and $E$.
Let $\triangle$ be the root system of $\g$, $Q$ be the root lattice and  
$$\g=\h\oplus (\bigoplus_{\alpha\in\triangle} \C e_\al)$$
be the weight decomposition. For any $\al\in Q$ we have a vertex operator
$$\Gamma_\al(z)=\sum_{n\in\Z} \Gamma_\al(n) z^{-n-1}$$
acting on the space $V_Q$. We let $V(\g)$ to denote the vertex operator algebra
associated with $\gh$.

\begin{theo}
\label{FK}
The identification $\Gamma_\al(n)\mapsto e_\al\T t^n$ defines an isomorphism
$V(\g)\simeq V_Q$, which sends a highest weight vector of $L_k$ to $|0\ket$.
\end{theo}

We finish this section with the description of principal subspaces 
of vertex operator algebras (see \cite{FS}).
Let $\al_1,\ldots,\al_N$ be a set of
of linearly independent vectors generating 
the lattice $Q$. Let $M=(m_{i,j})_{1\le i,j\le N}$ be the non degenerate 
matrix of the scalar
products of $\al_i$ ($m_{i,j}=(\al_i,\al_j)$) such that  
$m_{i,i}=2$ for all $i$. 
Consider the principal subspace   
$W_Q\hk V_Q$ generated from the vector
$|0\ket$ with an action of operators $\Gamma_{\al_i}(-n_i)$ with
$n_i\ge 1$  ($1\le i\le N$). Note that because of 
($\ref{htoh}$) one has $\Gamma_{\al_i}(-n) |0\ket =0$ for 
$n<1$.

Our goal is to describe $W_Q$ (in particular to find its character).  
We first realize this subspace as a quotient of a polynomial
algebra. Namely define $W'_Q$ as a quotient of 
$\C[a_i(-n)]$ with $1\le i\le N$, $n\ge 1$,  
by the ideal of relations generated with
$$a_i(z)^{(k)}a_j(z)^{(l)},\ k+l<m_{ij},$$
where $a_i(z)=\sum_{n\ge 1} z^n a_i(-n)$.
We note that 
$W_Q'=\bigoplus_{\bn\in\Z_{\ge 0}^N} W_{Q,\bn}'$, 
where $W_{Q,\bn}'$ is a subspace spanned by monomials in $a_i(k)$ such that
the number of factors of the type $a_{i_0}(k)$ with fixed $i_0$ is exactly
$n_{i_0}$. The $q$-character of $W_{Q,\bn}'$ is naturally defined by
$\deg_q a_i(k)=-k.$

The following lemma is standard (see for example \cite{FF}).
\begin{lem}
\label{charlem}
\begin{equation}
\label{charact}
\ch_q W_{Q,\bn}'=
\frac{q^{\bn M\bn/2}}{(q)_\bn},
\end{equation}
where $(q)_\bn=\prod_{j=1}^N (q)_{n_j}$, $(q)_n=\prod_{j=1}^n (1-q^j)$.
\end{lem}

The next proposition states  that the spaces $W_Q$ and $W_Q'$ are
isomorphic (see for example \cite{FF}).

\begin{prop}
\label{charvo}
The map $|0\ket\mapsto 1$, $\Gamma_{\al_i}(n)\mapsto a_i(n)$ 
induces the isomorphism
$$W_Q\simeq W_Q'.$$ 
In particular for any $\bn=(n_1,\ldots,n_N)\in\Z_{\ge 0}^N$
\begin{equation}
\label{chvo}
\ch_q (W_Q\cap \pi_{n_1\al_1+\ldots +n_N\al_N})=
\frac{q^{\bn M\bn/2}}{(q)_\bn}. 
\end{equation}
\end{prop}

\section{Level $1$ case}
\subsection{Algebras $A_1$ and $B_1$.}
We start with the description of the 
adjoint graded space $L_1^{gr}$ for vacuum
level $1$ $\slth$-module. Note that this case was considered in 
\cite{FFJMT} by the different
method. We make a connection to their approach in the end of the section.

The space  $L_1^{gr}$ carries the natural structure of the representation of abelian Lie
algebra $\slt^{ab}\T t^{-1}\C[t^{-1}]$ (see Remark $\ref{nab}$), 
where $\slt^{ab}$ is an abelinization of $\slt$, i.e.
the $3$-dimensional abelian Lie algebra. For $x\in\slt$ we denote the corresponding element of
$\slt^{ab}$ by $\tilde x$. We also set
$$
\tilde x(z)=\sum_{i<0} \tilde x_i z^{-i-1}.
$$
The space $L_1^{gr}$ is isomorphic to a quotient of the polynomial algebra 
$$\C[\tilde e_i, \tilde h_i,\tilde f_i]_{i<0}$$ by some ideal. 
Our goal is to show that the ideal of relations is generated with 
coefficients of the following series
\beq
\label{A1}
\tilde e(z)^2,\ \tilde e(z)\tilde h(z),\ 2\tilde e(z)\tilde f(z) -\tilde h(z)^2,\
\tilde h(z) \tilde f(z),\ \tilde f(z)^2, 
\eeq
i.e. by the $\slt$ consequences of the relation $\tilde e(z)^2$
($\slt$ acts on $\slt^{ab}\T t^{-1}\C[t^{-1}]$ via the adjoint action on the space
$\slt^{ab}$).
 
\begin{defin}
Let $A_1$ be an algebra generated with the commuting variables 
$\tilde x_i$, $x=e,h,f$, $i<0$ 
modulo the relations $(\ref{A1})$.
\end{defin}  

In order to make a connection between $A_1$ and $L_1^{gr}$  we need a modification of $A_1$. 

\begin{defin}
Let $B_1$ be an algebra generated with the abelian variables 
$\tilde x_i$, $x=e,h,f$, $i<0$ modulo the relations 
\beq
\label{B_1}
\tilde e(z)^2,\ \tilde e(z)\tilde h(z),\ \tilde h(z)^2,\
\tilde h(z)\tilde f(z),\ \tilde f(z)^2. 
\eeq
\end{defin}

We define the $(q,z,u)$-characters of $A_1$ and $B_1$ assigning the 
$(q,z,u)$-degree to each generator $\tilde x_i$:
$$
\deg_q \tilde x_i=-i,\ \deg_z \tilde e_i=2,\ \deg_z \tilde f_i=-2,\ \deg_z\tilde h_i=0,\
\deg_u \tilde x_i=1. 
$$

\begin{lem}
\label{LAB}
$\ch_{q,z,u} L_1^{gr}\le \ch_{q,z,u} A_1\le \ch_{q,z,u} B_1$
(the inequalities are true in each weight component).
\end{lem}
\begin{proof}
To show that $\ch_{q,z,u} L_1^{gr}\le \ch_{q,z,u} A_1$ it suffices to verify that relations
$(\ref{A1})$ hold in $L_1^{gr}$. In fact, $e(z)^2=0$ in $L_1$
and $\slt=\slt\T 1\hk \slth$ is acting on $L_1^{gr}$ as well as on $L_1$. 
In addition
all of the relations $(\ref{A1})$ can be produced from $\tilde e(z)^2$ by applying the operator
$f$.

To prove that $\ch_{q,z,u} A_1\le \ch_{q,z,u} B_1$  we introduce a filtration $G_s$ on $A_1$ by
setting $G_0$ to be the subspace generated with variables 
$\tilde e_i$, $\tilde f_i$ and
$$
G_{s+1}=\mathrm{span}\{\tilde h_i w:\ i<0, w\in G_s\}.
$$
Then the relation $2\tilde e(z)\tilde f(z) -\tilde h(z)^2=0$ (which holds in $A_1$) 
gives $\tilde h(z)^2=0$ in the adjoint graded space. Lemma is proved. 
\end{proof}
In the following subsections we give two proofs of the equality
$$\ch_{q,z,u} L_1 = \ch_{q,z,u} B_1.$$

\subsection{A basis of $B_1$.}
Recall the $ehf$-basis of $L_1$ (see Theorem $\ref{basis}$): 
\beq
\label{ehf}
 \dots e_{-2}^{c_2} f_{-1}^{a_1} h_{-1}^{b_1}e_{-1}^{c_1}
\eeq 
with 
\begin{enumerate}
\item[$(a)$] $a_i+a_{i+1}+b_{i+1}\le 1$,
\item[$(b)$] $a_i+b_{i+1}+c_{i+1}\le 1$,
\item[$(c)$] $a_i+b_i+c_{i+1}\le 1,$
\item[$(d)$] $b_i+c_i+c_{i+1}\le 1.$
\end{enumerate}
We now consider a modified set of restrictions
\begin{enumerate}
\item[$(a')$] 
$a_i+a_{i+1}+b_{i+1}\le 1,$
\item[$(b')$] 
$a_i+b_i+b_{i+1}\le 1,$
\item[$(c')$]
$b_i+b_{i+1}+c_{i+1}\le 1,$
\item[$(d')$]
$b_i+c_i+c_{i+1}\le 1,$
\end{enumerate}
and add one additional restriction
\begin{enumerate}
\item[$(N)$]
$b_i+a_{i+1}+c_{i+2}\le 2.$
\end{enumerate}
We refer to the monomials $(\ref{ehf})$ with restrictions $(a')$, $(b')$,
$(c')$, $(d')$ and $(N)$ as $ehf'$-monomials.

Conditions $(a')-(d')$ and $(N)$  can be expressed as follows
(the first picture explains $(a')-(d')$ and the second explains $(N)$):

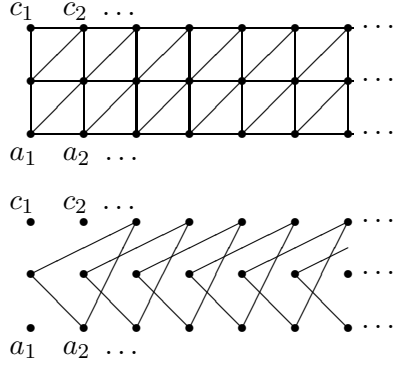
\begin{center}
\begin{picture}(180,70)

\multiput(40,20)(0,20){3}{\line(1,0){122}}
\multiput(40,20)(20,0){7}{\line(0,1){40}}

\multiput(40,20)(20,0){6}{\line(1,1){20}}
\multiput(40,40)(20,0){6}{\line(1,1){20}}
\multiput(165,20)(0,20){3}{\dots}
\multiput(40,20)(20,0){7}{\circle*{3}}
\multiput(40,40)(20,0){7}{\circle*{3}}
\multiput(40,60)(20,0){7}{\circle*{3}}
\put(32,10){$a_1$} 
\put(52,10){$a_2\ \dots$}
\put(32,65){$c_1$}
\put(52,65){$c_2 \ \dots$}

\end{picture}
\end{center}              

\begin{center}
\begin{picture}(180,70)


\multiput(40,40)(20,0){5}{\line(2,1){40}}
\multiput(140,40)(20,0){1}{\line(2,1){20}}
\multiput(40,40)(20,0){6}{\line(1,-1){20}}

\multiput(60,20)(20,0){5}{\line(1,2){20}}
\multiput(165,20)(0,20){3}{\dots}
\multiput(40,20)(20,0){7}{\circle*{3}}
\multiput(40,40)(20,0){7}{\circle*{3}}
\multiput(40,60)(20,0){7}{\circle*{3}}
\put(32,10){$a_1$} 
\put(52,10){$a_2\ \dots$}
\put(32,65){$c_1$}
\put(52,65){$c_2 \ \dots$}

\end{picture}
\end{center}              

\noindent The $ehf'$-monomials are monomials $(\ref{ehf})$
such that the sum of exponents over any
triangle in the first picture 
is less than or equal to $1$ and 
the sum of exponents over any
triangle in the second picture 
is less than or equal to $2$. 

\begin{lem}
\label{1}
The characters of $ehf$- and $ehf'$-monomials coincide,
\end{lem}
\begin{proof}
In order to prove our Lemma we construct a $(q,z,u)$-degree preserving 
bijection $\phi$ from the set
of $ehf$-monomials to the set of $ehf'$-monomials. Fix some $ehf$-monomial 
$m$ of the form $(\ref{ehf})$. If $m$ satisfies conditions 
$(a')$, $(b')$, $(c')$ and $(d')$ 
then we set $\phi(m)=m$. Now suppose that $j$ is a smallest number such that  
$(b')$ or $(c')$ is violated, i.e.
$$
a_j+b_j+b_{j+1}= 2 \text{ or } b_j+b_{j+1}+c_{j+1}= 2. 
$$
This means 
$$b_j=b_{j+1}=1,\ a_j=c_{j+1}=0.$$
We construct a new monomial 
$$m'=\dots f_{-n}^{a'_n} h_{-n}^{b'_n} e_{-n}^{c'_n}\dots  
f_{-1}^{a'_1} h_{-1}^{b'_1} e_{-1}^{c'_1},$$ 
which differs from $m$ only in terms $a_j, b_j, b_{j+1}, c_{j+1}$. 
The new values are given by
\beq
\label{transform}
a'_j=c'_{j+1}=1,\ b'_j=b'_{j+1}=0.
\eeq
The new monomial $m'$ satisfies conditions $(a')$-$(d')$ for all $i\le j$.
In addition it satisfies the condition $(N)$ for all $i<j$. In fact the 
violation of $(N)$ means that for some $i<j$ 
$b'_i=a'_{i+1}=c'_{i+2}=1$. Therefore $i=j-1$ and $b_{j-1}=b_j=1$. This
gives the violation of the conditions $(b')$ and $(c')$ for $i=j-1$, which
contradicts with the choice of $j$.

We repeat the procedure until all the conditions 
$(a')-(d')$ are satisfied for all $i$. Denote the result by $\phi(m)$. 
We note that by the same reason as above the condition $(N)$ is 
satisfied by $\phi(m)$. Therefore $\phi(m)$ is a $ehf'$-monomial.

We now show that $\phi$ is a bijection.
The inverse map $\phi^{-1}$ can be constructed in the following
way. Fix some $m'$, which is $ehf'$- but not $ehf$-monomial. Then for some $j$
$$
a'_j+b'_j+c'_{j+1}=2 \text{ or } a'_j+b'_{j+1}+c'_{j+1}=2,
$$  
which means $a'_j=c'_{j+1}=1$, $b'_j=b'_{j+1}=0$.
Then the inverse transformation to $(\ref{transform})$ is given by
$$
a_j=c_{j+1}=0,\  b_j=b_{j+1}=1.
$$
Lemma is proved.
\end{proof}

\begin{lem}
\label{2}
The set of $ehf'$-monomials spans $B_1$.
\end{lem}
\begin{proof}
For any ordered monomial $m$ of the form
\beq
\label{monom}
\dots f_{-n}^{a_n} h_{-n}^{b_n} e_{-n}^{c_n}\dots  
f_{-1}^{a_1} h_{-1}^{b_1} e_{-1}^{c_1} 
\eeq 
we denote by $\deg m$ the sum of all the exponents $a_i$, $b_i$ and $c_i$. 
Following \cite{FKLMM} we define a 
complete lexicographic ordering on the set of ordered monomials by the following rule.
If $\deg m>\deg m'$ then $m>m'$. Suppose $\deg m=\deg m'$. 
Then if $c_1<c'_1$ then
$m>m'$. If $\deg m=\deg m'$, $c_1=c'_1$ and $b_1< b'_1$ then $m>m'$. Next we compare $a_1$ 
and $a'_1$ and so on. 

We now show that any monomial $m$ which violates one of the conditions
$(a')$-$(d')$, $(N)$ can be rewritten as a sum of smaller monomials.
Suppose that for some $i$ the condition $(a')$ is violated, 
i.e. $a_i+a_{i+1}+b_{i+1}>1$. This means that at least two of numbers 
among $a_i, a_{i+1}, b_{i+1}$ are greater than $0$. But each of the products
$a_i a_{i+1}$, $a_i b_{i+1}$, $a_{i+1}b_{i+1}$ can be written as
a linear combination of smaller monomials using the relations
$$\sum_{\al+\be=2i+1} \tilde f_\al \tilde f_\be=0,\  
\sum_{\al+\be=2i+1} \tilde f_\al \tilde h_\be=0,\ 
\sum_{\al+\be=2i+2} \tilde f_\al \tilde h_\be=0$$
respectively.
By the similar reason the violation of  $(b')$, $(c')$ and $(d')$
allows to rewrite the corresponding monomial in terms of the smaller ones.
So we  finish the proof with rewriting a monomial which satisfies 
$(a')-(d')$ but not $(N)$. This reduces to rewriting a monomial
$\tilde h_i \tilde e_{i+1}\tilde f_{i+2}$ as a linear combination of a
smaller ones. We have
\begin{equation}
\label{rewr}
\tilde h_i \tilde e_{i+1}\tilde f_{i+2}=
-\tilde e_i \tilde h_{i+1} \tilde f_{i+2}+\dots=
\tilde e_i \tilde h_{i+2} \tilde f_{i+1}+\dots,
\end{equation}
where we use $\dots$ to denote a linear combination of the monomials, which
are smaller than $\tilde h_i \tilde e_{i+1}\tilde f_{i+2}$
(the relations $\tilde e(z)\tilde h(z)=0$ and $\tilde h(z)\tilde f(z)=0$
are used in $(\ref{rewr})$).
Now it is enough to note that the last expression is smaller than 
$\tilde h_i \tilde e_{i+1}\tilde f_{i+2}$. 
\end{proof}

\begin{cor}
\label{LB}
$\ch_{q,z,u} B_1=\ch_{q,z,u} L_1^{gr}$.
\end{cor}
\begin{proof}
From Lemmas $\ref{1}$ and $\ref{2}$ we know that the character of the set of 
$ehf'$-monomials is 
greater  than or equal to the character of
$B_1$ and coincides with the character of $L_1$. Now our Corollary 
follows from Lemma  $\ref{LAB}$. 
\end{proof}

\subsection{The degeneration procedure.}
In this section we give an alternative proof of Corollary $\ref{LB}$.
We involve the degeneration procedure which will be used later in the case
of general level $k$. 

Let $p,q,r\in \R^3$ be linearly independent vectors with the scalar
products 
$$(p,p)=(q,q)=(r,r)=2, (p,q)=(q,r)=1, (p,r)=0.$$
Let $Q_1$ be a lattice generated by $p,q,r$ and $\Gamma_p(z)$, $\Gamma_q(z)$, 
$\Gamma_r(z)$ be the corresponding vertex operators.

\begin{lem}
\label{voB1}
The identification 
$$\Gamma_p(n)\to \tilde e_n,\ 
\Gamma_q(n)\to \tilde h_n,\ \Gamma_r(n)\to \tilde f_n$$
provides an isomorphism $W_{Q_1}\simeq B_1$, $|0\ket\to 1$. 
\end{lem}
\begin{proof}
Follows from the Proposition $\ref{charvo}$.
\end{proof}

We now consider the Lie algebra $\widehat{\s_4}$ and its vertex operator 
realization
provided by the Frenkel-Kac construction.  For $1\le i,j\le 4$ we set
$$E_{i,j}(z)=\sum_{n\in\Z} (E_{i,j}\T t^n)z^{-n-1}.$$
Using Theorem $\ref{FK}$ and Lemma $\ref{voB1}$ we obtain the following Lemma. 

\begin{lem}
The identification 
$$\tilde e_n\to E_{2,3}\T t^n,\ 
\tilde h_n\to E_{2,4}\T t^n,\ \tilde f_n\to E_{1,4}\T t^n$$
provides an isomorphism between $B_1$ and the subspace of $V(\s_4)$ generated
from the highest weight vector with the fields $E_{2,3}(z)$, $E_{2,4}(z)$ and
$E_{1,4}(z)$. 
\end{lem}

We want to show that the character of $B_1$ is smaller than or equal to
the character of $L_1$. In order to do this we construct a deformation
of the Lie algebra $\slt\hk \s_4$ to the subspace spanned by
$E_{2,3}$, $E_{2,4}$ and $E_{1,4}$.

\begin{prop}
\label{degen}
There exists a continuous family of Lie subalgebras 
$S(\veps)$ of $\s_4$, $0\le\veps\le 1$ such that
\begin{enumerate}
\item[$a).$] $S(\veps)\simeq \slt$ for $0\le \veps< 1$,
\item[$b).$] $S(1)$ is spanned with $E_{2,3}$, $E_{2,4}$ and $E_{1,4}$,
\item[$c).$] There exist standard basises $e(\veps)$, $h(\veps)$, $f(\veps)$ of
$S(\veps)$ ($0\le\veps< 1$) such that 
$$\lim_{\veps\to 1} \C e(\veps) =\C E_{1,4},\ 
\lim_{\veps\to 1} \C h(\veps) =\C E_{2,4},\
\lim_{\veps\to 1} \C f(\veps) =\C E_{2,3}.$$
\end{enumerate}
\end{prop}
\begin{proof}
Let $v_1, v_2, v_3, v_4$ be a basis of $\C^4$. For $0\le \veps <1$ we let 
$S(\veps)$ to denote
the subalgebra determined by the following conditions:
\begin{itemize}
\item $S(\veps)$ preserves $\mathrm{span}\{v_1,v_2\}$,
\item 
$S(\veps)$ annihilates 
$\mathrm{span}\{(1-\veps)^{\frac{1}{2}}v_3+ \veps v_1, 
(1-\veps)v_4+\veps v_2\}$,
\item 
$S(\veps)$ contains only traceless matrices.
\end{itemize}
Then $S(\veps)$ consists of the matrixes of the form
$$
\begin{pmatrix}
x & y & \frac{-x\veps}{\sqrt{1-\veps}} & \frac{-y\veps}{1-\veps}\\
z & -x & \frac{-z\veps}{\sqrt{1-\veps}} & \frac{x\veps}{1-\veps}\\
0 & 0 & 0 & 0 \\
0 & 0 & 0 & 0 
\end{pmatrix}.
$$ 
Let $e(\veps)$, $h(\veps)$, $f(\veps)$ be the standard basis of 
$S(\veps)\simeq \slt$ (we fix the identification via the upper left $2\times 2$
corner of $S(\veps)$). Then
$$\lim_{\veps\to 1} (\veps-1) e(\veps) = E_{1,4},\ 
\lim_{\veps\to 1} (1-\veps) h(\veps) = E_{2,4},\
\lim_{\veps\to 1} (1-\veps)^{\frac{1}{2}} f(\veps) = -E_{2,3}.
$$
This finishes the proof of the Proposition.
\end{proof}

\begin{cor}
\label{L1toB1}
There exists a continuous family of subspaces $M(\veps)\hk V(\s_4)$, 
$0\le\veps\le 1$ such that
$M(\veps)\simeq L_1$ for $0\le\veps< 1$ and $M(1)\simeq B_1$. 
\end{cor}

\begin{cor}
\label{LA}
$L_1^{gr}$ is isomorphic to $A_1$ as a representation of 
the algebra 
$\slt^{ab}\T t^{-1}\C[t^{-1}]$.
\end{cor}
\begin{proof}
Follows from Lemma $\ref{LAB}$ and Lemma $\ref{LB}$.
\end{proof}

Recall the $(q,z,u)$-character of $L_{0,1}^{gr}$ by 
\beq
\ch_{q,z,u} L_1^{gr}=\sum_{s=0}^\infty u^s \ch_{q,z} (F_s/F_{s-1}). 
\eeq

\begin{prop}
\beq
\label{ferm}
\ch_{q,z,u} L_1^{gr}=
\sum_{n^+,n^0,n^-\ge 0} u^{n^++n^0+n^-} z^{2(n^+ -n^-)}
\frac{q^{(n^+)^2+(n^0)^2+(n^-)^2+n^+n^0+n^0n^-}}
{(q)_{n^+}(q)_{n^0}(q)_{n^-}},
\eeq
where $(q)_n=\prod_{i=1}^n (1-q^i)$.
\end{prop}

\begin{proof}
We first calculate the character of $B_1$.
Consider the map $\phi$ from $B_1^*$ to the space of polynomials in $3$ groups of variables:
\begin{multline*}
(\phi(\theta))(x^-_1,\dots ,x^-_{n^-}, x^0_1,\dots ,x^0_{n^0}, x^+_1,\dots, x^+_{n^+})=\\
\theta(\tilde f(x^-_1)\dots \tilde f(x^-_{n^-}) \tilde h(x^0_1)\dots \tilde h(x^0_{n^0}) 
\tilde e(x^+_1)\dots \tilde e(x^+_{n^+})).
\end{multline*}
Because of the relations $(\ref{B_1})$ the image of $\phi$  coincides with the 
the space of
polynomials of the form
\begin{multline}
\label{dualsp}
\prod_{\at{\alpha=+,-,0}{1\le i^s\le n^\alpha}} x^\alpha_i 
\prod_{\at{\alpha=+,-,0}{1\le i< j\le n^\alpha}} (x^\alpha_i-x^\alpha_j)^2
\prod_{\at{1\le i\le n^+}{1\le  j\le n^0}} (x^+_i-x^0_j)
\prod_{\at{1\le i\le n^0}{1\le  j\le n^-}} (x^0_i-x^-_j)\times\\ 
F(x^-_1,\dots ,x^-_{n^-}, x^0_1, \dots ,x^0_{n^0}, x^+_1, \dots,  x^+_{n^+}),
\end{multline}
where $F$ is an arbitrary polynomial in $x^\alpha_i$, symmetric in each group ($+,0,-$) of variables. 
The natural $(q,z,u)$-grading on the space $B_1^*$ comes from the grading on $B_1$.
It is easy to see that the corresponding character of the space of polynomials $(\ref{dualsp})$
is given by
$$
u^{n^++n^0+n^-} z^{2(n^+ -n^-)}
\frac{q^{(n^+)^2+(n^0)^2+(n^-)^2+n^+n^0+n^0n^-}}
{(q)_{n^+}(q)_{n^0}(q)_{n^-}}.
$$
Now our Lemma follows from Corollary $\ref{LB}$.
\end{proof}

\begin{cor}
\beq
\ch_{q,z} L_1=
\sum_{n^+,n^0,n^-\ge 0} z^{2(n^+ -n^-)}
\frac{q^{(n^+)^2+(n^0)^2+(n^-)^2+n^+n^0+n^0n^-}}
{(q)_{n^+}(q)_{n^0}(q)_{n^-}}.
\eeq
In particular the right hand side equals to the well known expression 
for the character of $L_1$:
$$
\ch_{q,z} L_1=\frac{1}{(q)_{\infty}}\sum_{n\in\Z} z^{2n} q^{n^2},
$$
where $(q)_\infty=\prod_{i\ge 1} (1-q^i)$
\end{cor}

We now compare our formula with one from \cite{FFJMT}, which is given in terms of
supernomial coefficients (see \cite{SW}). 

\begin{cor}
\beq
\label{fus}
\ch_{q,z,u} L_1^{gr}=
\sum_{m\ge 0} \frac{u^m q^{m^2}}{(q)_m} \sum_{-m\le l\le m} 
z^{2l} S_{m,l}(q), 
\eeq
where 
$$S_{m,l}(q)=\sum_{m\ge \nu\ge m-l-\nu} q^{(\nu+l-m)(\nu+l)+\nu(\nu-m)} 
\qb{m}{\nu} \qb{\nu}{m-l-\nu}
$$
and $\qb{n}{m}=\frac{(q)_n}{(q)_m(q)_{n-m}}$. 
\end{cor}
\begin{proof}
We note that
$$
\frac{1}{(q)_m} \qb{m}{\nu} \qb{\nu}{m-l-\nu}=\frac{1}{(q)_{m-\nu}(q)_{2\nu-m+l}(q)_{m-l-\nu}}.
$$
The change of variables
$$m-\nu=n_+,\ 2\nu-m+l=n_0,\ m-l-\nu=n_-$$
identifies formulas $(\ref{fus})$ and $(\ref{ferm})$.
\end{proof}

\section{The general case}
In this section 
we consider the PBW-filtration on the vacuum level $k$ representation 
$L_k$ and 
the corresponding adjoint graded space $L_k^{gr}$. 

\subsection{Algebras $A_k$ and $B_k$.}
We introduce $2$ series of algebras $A_k$ and $B_k$
generated with Fourier coefficients of the abelian currents $\tilde e(z)$, $\tilde f(z)$ and 
$\tilde h(z)$, where 
$$\tilde x(z)=\sum_{i\ge 1} z^{-i-1} \tilde x_i.$$

\begin{defin}
An algebra $A_k$ is a quotient of the polynomial algebra in variables 
$\tilde e_i$, $\tilde h_i$, $\tilde f_i$, $i<0$,
by the ideal generated by Fourier coefficients of $\tilde e(z)^{k+1}$ and all its $\slt$ 
consequences, i.e. by coefficients of  $2k+3$ series: 
\beq
\label{A_krel}
\tilde e(z)^{k+1}, \ \tilde e(z)^k\tilde h(z),\ 
\tilde e(z)^{k-1}\tilde h^2(z)-2\tilde e(z)^k\tilde f(z), \dots, \tilde f(z)^{k+1}. 
\eeq   
\end{defin}
In other words these relations can be described as follows. Identify $\slt$ with its $3$-dimensional
irreducible representation $\pi_2$ with $e$ being the highest weight vector. Then we have an embedding
$$\imath:\pi_{2k+2} \hk\pi_2^{\T (k+1)}$$
of $(2k+3)$-dimensional irreducible $\slt$ module $\pi_{2k+2}$ into the tensor power 
$\pi_2^{\T (k+1)}$
(the image of $\pi_{2k+2}$ is generated from $e^{\T (k+1)}\in  \pi_2^{\T (k+1)}$ by the action
of universal enveloping algebra of $\slt$). Define an affinization map $\alpha$, which sends an
element of the tensor power $\pi_2^{\T (k+1)}$ to the product of the corresponding series: 
$$
\alpha (x^1\T\dots\T x^{k+1})= \tilde x^1(z)\dots \tilde x^{k+1}(z),\ x^i= e,h,f. 
$$
Then the defining relations of $A_k$ are coefficients of  $\alpha (\imath ( \pi_{2k+2}))$.

We note that $(q,z,u)$-character of $A_k$ is naturally defined by
\beq
\label{deg}
\deg_u \tilde x_i=1,\  \deg_q \tilde x_i=-i,\ 
\deg_z \tilde e_i=2,\ \deg_z \tilde h_i=0,\deg_z \tilde f_i=-2.
\eeq

\begin{lem}
\label{LkAk}
$\ch_{q,z,u} L_k^{gr}\le \ch_{q,z,u} A_k.$
\end{lem}
\begin{proof}
Follows from the equality $e(z)^{k+1}=0$ in $L_k$.
\end{proof}

\begin{defin}
An algebra $B_k$ is a quotient of 
the polynomial algebra in variables $\tilde e_i$, $\tilde h_i$, $\tilde f_i$, $i\le -1$
by the ideal generated by Fourier coefficients of $2k+3$ series
\beq
\label{B_k}
\tilde e(z)^i \tilde h(z)^{k+1-i},\ 1\le i\le k+1;\ 
\tilde h(z)^i \tilde f(z)^{k+1-i},\ 0\le i\le k+1.
\eeq
\end{defin}

We note that $(q,z,u)$-character of $B_k$ is naturally defined by $(\ref{deg})$.
\begin{lem}
\label{AkBk}
$\ch_{q,z,u} A_k\le \ch_{q,z,u} B_k.$
\end{lem}

\begin{proof}
We introduce a filtration $G_s$ on $A_k$ by
setting $G_0$ to be the subspace generated by variables $\tilde e_i$, $\tilde f_i$ (but not
$\tilde h_i$) and
$$
G_{s+1}=\mathrm{span}\{\tilde h_i w:\ i<0, w\in G_s\}.
$$
Then for $0\le i\le k+1$ the relation 
$$\alpha (\imath (f^i \cdot e^{\T (k+1)}))=0$$ 
which holds in $A_k$ contains a term $\tilde e(z)^{k+1-i}\tilde h(z)^i$ and the coefficients
of the difference
$$\alpha(\imath (f^i \cdot e^{\T (k+1)}))-\tilde e(z)^{k+1-i}\tilde h(z)^i$$ 
belongs to $G_{i-1}$. This means that the relation $ \tilde e(z)^{k+1-i}\tilde h(z)^i=0$ 
holds in the adjoint graded space
$G_0\oplus\bigoplus_{s>0} (G_s/G_{s-1})$. Similarly the rest of the relations $(\ref{B_k})$ are true
in the adjoint graded space. Lemma is proved.
\end{proof}

\begin{cor}
\label{LkBk}
$\ch_{q,z,u} L_k^{gr}\le \ch_{q,z,u} B_k.$
\end{cor}

\subsection{A quadratic algebra $C_k$ and principal subspace $D_k$.}
We consider a set of commuting variables 
$$\tilde x^{[l]}_i,\  x= e, h, f, \ 1\le l\le k,\  i\le -l$$ 
and the corresponding currents
$\tilde x^{[l]}(z)=\sum_{i< 0} z^{-i-l} \tilde x^{[l]}_i$ 
(we set $\tilde x^{[l]}_i=0$ for $i>-l$).
For the series $p(z)$ let $p(z)^{(r)}$ be the $r$-th derivative.

\begin{defin}
Let $C_k$ be the quotient of the polynomial algebra in commuting variables
$\tilde x^{[l]}_i$, $1\le l\le k$, $i\le -l$ by the ideal of relations generated by coefficients of 
currents
\begin{gather}
\label{min}
\tilde x^{[l]}(z)^{(\alpha)} \tilde x^{[m]}(z)^{(\beta)} \text{ for } x= e, h, f,\ 
\alpha+\beta<2\min(l,m),\\
\label{max1}
\tilde e^{[l]}(z)^{(\alpha)} \tilde h^{[m]}(z)^{(\beta)} \text{ for } 
\alpha+\beta<\max(0, l+m-k),\\
\label{max2}
\tilde h^{[l]}(z)^{(\alpha)} \tilde f^{[m]}(z)^{(\beta)} \text{ for } 
\alpha+\beta<\max(0,l+m-k).
\end{gather}  
\end{defin}

We define the $(q, z, u)$-degree on $C_k$ by the formulas
$$
\deg_q \tilde x^{[l]}_i=-i,\ \deg_z \tilde e^{[l]}_i=2l,\ \deg_z \tilde h^{[l]}_i=0,\ 
\deg_z \tilde f^{[l]}_i=-2l,\ \deg_u \tilde x^{[l]}_i=l.
$$
In what follows we show that the $(q,z,u)$-characters of $B_k$ and $C_k$ coincide. 

We will need the following Lemma from \cite{FS, FJKLM}.

\begin{lem}
\label{filtr}
Consider the quotient of the polynomial algebra
$\C[a_0,a_{-1},\ld]$ by the ideal generated with the coefficients of the series 
$a(z)^{k+1}$. Then there exists a filtration
$F_\mu$ of this quotient (labeled by Young diagrams $\mu$)
such that the adjoint graded algebra is
generated with the coefficients of series $a^{[i]}(z)$, which are the images of
powers $a(z)^i$, $1\le i\le k$. 
In addition defining relations in the adjoint graded algebra
are of the form
\begin{equation}
\label{rel}
a^{[i]}(z)^{(l)} a^{[j]}(z)^{(r)}=0 \text{ if } l+r< 2\min(i,j).
\end{equation}
\end{lem}

\begin{lem}
\label{BC}
$\ch_{q,z,u} B_k\le \ch_{q,z,u} C_k$.
\end{lem}
\begin{proof}
Using Lemma $\ref{filtr}$ we define a filtration on $B_k$ such that the adjoint graded space is 
generated
by the images of coefficients of the series $\tilde x(z)^l$, $x=e, h, f$, 
$1\le l\le k$. Denote the corresponding series by $\tilde x^{[l]}(z)$ and the corresponding Fourier
coefficients by $\tilde x^{[l]}_i$. Then from Lemma $\ref{filtr}$ we obtain the relations 
$(\ref{min})$. The
relations $(\ref{max1})$, $(\ref{max2})$  in the adjoint graded space follows from the relations 
\begin{gather}
(\tilde e(z)^l)^{(\alpha)} (\tilde h(z)^m)^{(\beta)}=0 \text{ for } 
\alpha+\beta<\max(0, l+m-k),\\
(\tilde h(z)^l)^{(\alpha)} (\tilde f(z)^m)^{(\beta)}=0 \text{ for } 
\alpha+\beta<\max(0,l+m-k),
\end{gather}  
which hold in $B_k$. We thus obtain that all of the relations of $C_k$ are true in the adjoint
graded space of $B_k$ with respect to the certain filtration. Lemma is proved.
\end{proof}  

We want to show that 
$\ch_{q,z,u} B_k\ge \ch_{q,z,u} C_k$.
We use the vertex operator technique. 
Fix an integer $N$ such that there exists a set of linearly independent vectors
$p_i, q_i, r_i\in\h=\R^N$, $1\le i\le k$ with the scalar products
\beq
\label{scal}
(p_i,p_j)=(q_i,q_j)=(r_i,r_j)=2\delta_{i,j},\ (p_i,q_j)=(q_i,r_j)=\delta_{i, k+1-j},\ (p_i,r_j)=0.
\eeq
For example, setting $N=3k$ and fixing some orthonormal basis $e_i$ with respect to 
$(\cdot,\cdot)$ one can define
$$
p_i=\sqrt{2}e_i ,\ q_i=\frac{1}{\sqrt{2}} e_{k+1-i} + \sqrt{\frac{3}{2}} e_{k+i},\
r_i=\sqrt{\frac{2}{3}} e_{2k+1-i}+ \sqrt{\frac{4}{3}} e_{2k+i}.
$$ 
We consider a lattice $Q$ generated by the vectors $p_i$, $q_i$, $r_i$ and the corresponding VOA
$V_Q$. Set
\beq
\label{V}
V_e(z)=\sum_{i=1}^k \Gamma_{p_i}(z),\ V_h(z)=\sum_{i=1}^k \Gamma_{q_i}(z),\
V_f(z)=\sum_{i=1}^k \Gamma_{r_i}(z).  
\eeq
Note that due to the Proposition $\ref{vorel}$ one has
$$
\Gamma_{p_i}(z)\Gamma_{q_{k+1-i}}(z)=0=\Gamma_{q_i}(z)\Gamma_{r_{k+1-i}}(z).
$$ 
Therefore for 
any $0\le i\le k+1$
\beq
\label{Dk}
V_e(z)^i V_h(z)^{k+1-i}=V_h(z)^iV_f(z)^{k+1-i}=0.
\eeq
We let $D_k$ to denote the space generated with the Fourier 
coefficients of 
$V_e(z)$, $V_h(z)$, $V_f(z)$ from the highest weight vector.

\begin{lem}
\label{BCD}
$\ch_{q,z,u} B_k\ge \ch_{q,z,u} D_k\ge \ch_{q,z,u} C_k$.
\end{lem} 
\begin{proof}
The equation $(\ref{Dk})$ proves the first part of our Lemma. To show the second
enequlity we use the degeneration procedure.
Namely we construct a deformation $D_k(\veps)$  of $D_k$ such that 
$D_k(\veps)\simeq D_k$ for 
$1\ge \veps > 0$ and $D_k(0)$ contains $C_k$. Namely let 
$D_k(\veps)$ be the subspace generated from the highest weight vector by the 
Fourier
coefficients of $V_e(\veps z)$, $V_h(\veps z)$, $V_f(\veps z)$. Then the limit
$\lim_{\veps\to 0} D_k(\veps)$ contains the Fourier coefficients of the series
\begin{gather*}
\Gamma_{p_1+\dots + p_l}(z)=\lim_{\veps\to 0} 
\veps^{\frac{l(l-1)}{2}} (V_e(\veps z))^l,\\
\Gamma_{q_1+\dots + q_l}(z)=\lim_{\veps\to 0} 
\veps^{\frac{l(l-1)}{2}} (V_h(\veps z))^l,\\
\Gamma_{r_1+\dots + r_l}(z)=\lim_{\veps\to 0} 
\veps^{\frac{l(l-1)}{2}} (V_f(\veps z))^l.
\end{gather*}
We note that 
\begin{gather*}
(p_1+\dots + p_i, p_1+\dots + p_j)= (q_1+\dots + q_i, q_1+\dots + q_j)=2\min(i,j),\\ 
(r_1+\dots + r_i, r_1+\dots + r_j)=2\min(i,j),\\
(p_1+\dots + p_i, q_1+\dots + q_j)= \max(0, i+j-k),\\
(r_1+\dots + r_i, q_1+\dots + q_j)= \max(0, i+j-k).
\end{gather*}
Let $Q_k$ be the lattice generated by the vectors
$$
p_1+\dots +p_l,\ q_1+\dots + q_l,\ r_1+\dots + r_l,\ 1\le l\le k.
$$
Then using Proposition $\ref{charvo}$ we obtain that the 
principal subspace $W_{Q_k}$ is isomorphic to $C_k$.
The isomorphism is given by the identification
$$
\Gamma_{p_1+\dots +p_l}(i)\mapsto \tilde e^{[l]}_i,\
\Gamma_{q_1+\dots +q_l}(i)\mapsto \tilde h^{[l]}_i,\
\Gamma_{r_1+\dots +r_l}(i)\mapsto \tilde f^{[l]}_i.
$$
This gives 
$$
\ch_{q,z,u} D_k\ge \ch_{q,z,u} C_k.
$$
Lemma is proved. 
\end{proof}

Lemmas $\ref{BC}$ and $\ref{BCD}$ give the following Corollary:
\begin{cor}
\label{B=C=D}
$\ch_{q,z,u} B_k= \ch_{q,z,u} C_k= \ch_{q,z,u} D_k$.
\end{cor}

\begin{prop}
\label{LkDk}
$\ch_{q,z,u} L_k\ge \ch_{q,z,u} D_k$. 
\end{prop}
\begin{proof}
We recall that there exists an embedding $L_k\hk L_1^{\T k}$ such that
the highest weight vector $v_k$ of $L_k$ maps to the tensor power $v_1^{\T k}$
and for any $x\in\slt$ the current $x(z)$ on $L_k$ corresponds to the sum
$\sum_{i=1}^k  x^{(i)}(z)$, where
$$x^{(i)}(z)=\Id\T\dots\T x(z)\T\dots\T \Id$$
($x(z)$ on the $i$-th place). From the definition of $D_k$ we also have an
embedding $D_k\hk D_1^{\T k}$ (see the definition $(\ref{V})$). Now the
degeneration from the Corollary $\ref{L1toB1}$  gives the degeneration of 
$L_k$. From the part $c)$ of the 
Proposition $\ref{degen}$ we conclude that the limit of this degeneration 
contains $D_k$. Proposition is proved.
\end{proof}

\begin{theo} $\phantom{a}$
\begin{itemize}
\item The $(q,z,u)$ characters of $L_k^{gr}$, $A_k$, $B_k$, $C_k$ and $D_k$ 
coincide.
\item $L_k^{gr}\simeq A_k$ as the modules over the abelian algebra with a basis 
$\tilde e_i$, $\tilde h_i$, $\tilde f_i$, $i\le -1$.
\end{itemize}
\end{theo}
\begin{proof}
Follows from Lemmas $\ref{LkAk}$, $\ref{AkBk}$, Corollary $\ref{B=C=D}$ 
and Proposition $\ref{LkBk}$.
\end{proof}

\subsection{The character formula.}
In this section we compute the $(q,z,u)$-character of $L_k^{gr}$ .

\begin{prop}
\begin{multline}
\label{charform}
\ch_{q,z,u} C_k=\sum_{\bn^-,\bn^0,\bn^+\in \Z_{\ge 0}^k}
u^{|\bn^+|+|\bn^0|+|\bn^-|} z^{2(|\bn^+|-|\bn^-|)}\times\\ 
\frac{ 
q^{\frac{1}{2}(\bn^+ A \bn^+ + \bn^0 A \bn^0 + \bn^- A \bn^-) + \bn^+ B \bn^0 + \bn^0 B \bn^-}}
{(q)_{\bn^+}(q)_{\bn^0}(q)_{\bn^-}},
\end{multline}
where for $\bn\in\Z_{\ge 0}^k$ we set $|\bn|=\sum_{i=1}^k i \bn_i$ and matrixes $A$ and $B$ are 
defined by 
$$A_{i,j}=2\min(i,j),\ B_{i,j}=\max(0, i+j-k).$$
\end{prop}

\begin{proof}
Follows from the vertex operator realization of $C_k$ constructed in 
Corollary $\ref{B=C=D}$
and Proposition $\ref{charvo}$.
\end{proof}

\begin{theo}
The $(q,z,u)$-character of $L_k^{gr}$ is given by the right hand side of
$(\ref{charform})$.
\end{theo}

\newcounter{a}
\setcounter{a}{1}

\end{document}